

\magnification=1200
\nopagenumbers
\parindent= 15pt
\baselineskip=13pt

\hsize=13.75cm
\vsize=19.59cm
\hoffset=-0.15cm

\input amssym.def
\input amssym.tex

 at6.5pt
\font\srm=cmr8

\font\csc=cmcsc10
\font\title=cmr14 at 16pt

\font\teneusm=eusm10    
\font\seveneusm=eusm7  
\font\fiveeusm=eusm5    
\newfam\eusmfam

\textfont\eusmfam=\teneusm 
\scriptfont\eusmfam=\seveneusm
\scriptscriptfont\eusmfam=\fiveeusm

\def\Re{{\rm Re}\,}
\def\Im{{\rm Im}\,}

\def\scr#1{{\scriptstyle{#1}}}
\def\c#1{{\cal #1}}

\def\B#1{{\Bbb #1}}

\def\rightheadline{\hfil{\srm Extension of
the Linnik Phenomenon}
\hfil\tenrm\folio}
\def\leftheadline{\tenrm\folio\hfil{\srm Y. Motohashi}\hfil}
\def\emptyheadline{}
\headline{\ifnum\pageno=1 \emptyheadline\else
\ifodd\pageno \rightheadline \else \leftheadline\fi\fi}

\def\firstpage{\hss{\vbox to 1.0cm
{\vfil\hbox{\rm\folio}}}\hss}
\def\emptyfootline{\hfil}
\footline{\ifnum\pageno=1\firstpage\else
\emptyfootline\fi}

\centerline{\title An Extension of the
Linnik Phenomenon} 
\vskip 0.7cm
\centerline{\csc Yoichi Motohashi}
\vskip 1cm 
\noindent
{\bf Abstract:} We shall prove an extension of Yu.V. 
Linnik's phenomenon concerning
C.L. Siegel's exceptional zeros of Dirichlet $L$-functions.
Also we shall prove a new
version of the Brun--Titchmarsh theorem. These two 
subjects are closely related to each other. 
The basic tool is our theory of A. Selberg's 
$\Lambda^2$-sieve, which is developed in 
our old Tata lecture notes [12].
Constants which are involved either explicitly 
or implicitly in our discussion are all 
effectively computable, save for the one
in the statement of a theorem of Siegel
which has no relation with our main assertions.
We add that the discussion in
Sections 1--4 is a substantially revised and 
augmented version of Section 9.3 of our recent
monograph [15, Vol.\ I] and that of
Section 5 an extraction
from our article [16] posted in arXiv.
\bigskip
\noindent
{\bf 1. Our motivation in a historical perspective.}  
We shall first make explicit our notion 
of exceptional zeros: Let $\chi$ denote 
a genetic Dirichlet character, 
with which the $L$-function 
$L(s,\chi)$ is associated. Let $T>0$ 
be sufficiently large, and let $Z_T=\{\rho=
\beta+i\gamma\}$ be the set of all non-trivial zeros $\rho$
in the region $|\Im s|\le T$ of the function
$\prod_{q\le T}\prod^*_{\chi\bmod q}L(s,\chi)$; 
here and in what follows the asterisk means 
that relevant characters are primitive. 
Then we have that there exists an
effectively computable absolute constant $a_0>0$
such that
$$
\hbox{$\displaystyle
\max_{\rho\in Z_T}\beta< 1-{a_0\over\log T},\;$
save for a possibly existing $\beta_T\in Z_T$.
}\eqno(1.1)
$$
If $\beta_T$ exists, it should be real and simple, 
and we designate both itself and
the relevant unique primitive character $\chi_T$, with
$L(\beta_T,\chi_T)=0$, 
as $T$-exceptional. It is known that $\chi_T$ is real.
As a matter of fact, we have more precisely
$$
1-{a_0\over \log T}\le\beta_T
\le 1-{1\over \sqrt{T}(\log T)^3}.\eqno(1.2)
$$
The upper bound is a  
consequence of the Dirichlet class number formula
for  quadratic number fields. However, the formula
can be dispensed with. See [7, Section 2, Chapter IX] as well as
[15, Section 4.5, Vol.\ I]; in both monographs are 
comprehensible accounts
of the theory of the zeta and Dirichlet $L$-functions and the theory of
the distribution of primes. 
\par
One of the
most tantalising problems in number theory is the elimination
of the possibility of the existence of exceptional zeros.
It is generally believed that they
do not exist at all. A way to confirm this 
is to improve the Brun--Titchmarsh theorem in the
manner to be made explicit at a later part of Section 5, 
which remains, however, to be one of the most difficult problems in
number theory. Also one may have a hope 
in employing our Theorem 1 for this purpose; we shall 
try to be precise in one of the remarks there.
\par
Exceptional zeros,
if exist, would mar the quality of the prime
number  theorem for arithmetic progressions,
especially  when the uniformity is
taken into account with respect to varying 
moduli. Whence is acclaimed Siegel's assertion (1936)
that for any fixed $\varepsilon>0$ there exists a 
$c(\varepsilon)>0$ such that
$$
\beta_T<1-{c(\varepsilon)\over T^\varepsilon},
\eqno(1.3)
$$
since this yields a fine uniformity as far as moduli
are relatively small yet meaningful as demonstrated
in a variety of fundamental works
such as I.M. Vinogradov's resolution (1937)
of the ternary
Goldbach conjecture and the mean prime number theorem
of E. Bombieri and A.I. Vinogradov (1965, independently).
\par
However, $c(\varepsilon)$ is ineffective; 
that is, all
known proofs of Siegel's theorem assert only 
the existence of $c(\varepsilon)$
for each $\varepsilon$ and do not provide 
any means to evaluate its actual values. 
This causes severe difficulties 
in various basic problems. The most outstanding
among them is the estimation of  the size 
of the least prime that appears in a 
given  arithmetic progression. 
To resolve this particular difficulty,  Linnik (1944)
greatly refined Siegel's theorem $(1.3)$ by providing
a  quantitative version of 
the Deuring--Heilbronn 
phenomenon or
the repelling effect of an exceptional zero
towards other zeros that had been 
observed in the efforts to solve C.F. Gauss' conjecture  
concerning the class
numbers of imaginary quadratic 
number fields. Linnik's
theorem  or rather his phenomenon asserts that the inequality
$$
\beta\le1-{a_0\over\log T}\log{a_0 e\over
(1-\beta_T)\log T}\eqno(1.4)
$$
holds for all $\rho\ne\beta_T$ in $Z_T$.
With this and a zero density result of a special type 
which is another important contribution by him, 
Linnik could reach his famed Least Prime Number 
Theorem for arithmetic progressions.
\par
Linnik's argument is, however, 
exceptionally involved, especially in his intriguing use of an
analytic convexity argument.
After K.A. Rodosskii's simplification, P. Tur\'an and S. Knapowski 
(1961-62) developed a relatively accessible alternative 
approach to Linnik's two assertions via Tur\'an's power sum method
which superseded the convexity argument.
Their theory was later improved considerably by E. Fogels (1965).
\par 
Then, a totally different approach was devised by A. Selberg  
(1973-74), which was in essence yet another replacement of the
convexity argument by an infusion 
of his quasi-characters into the theory of $L$-functions and the theory
of the large sieve of Linnik; see Th\'eor\`eme 7 of [1]. We 
observed later that Selberg's quasi-characters
come from the optimal $\lambda$-weights for the simplest
situation of his $\Lambda^2$-sieve, that is, sifting integers
in an interval with residue classes to be discarded being equal
to $0$. With this, we were able to
extend the notion of quasi-characters
by considering the $\Lambda^2$-sieve applied to sequences
of values of a class of multiplicative functions.
Combined the result with the large sieve, we
obtained a fairly simplified proof [11] of 
Linnik's least prime number theorem or more precisely
its extended version due to P.X. Gallagher (1970) who had
used the power sum method following Fogels. In passing,
we mention that M. Jutila [6] developed an argument
which is similar to but less general than ours.
\medskip
In retrospect, the r\^ole played by the sieve method in our
argument is in one way an
enhancement of a sieve aspect in Linnik's method as well as
Tur\'an--Knapowksi's, Fogels' and Gallagher's. Namely
those people applied the Brun--Titchmash theorem at respective
crucial stages in their arguments. 
\medskip
Now, the aim of the present work is to extend the Linnik 
phenomenon (1.4) to zeros of $L$-functions 
which may belong to a family of functions much wider than
hitherto considered, by employing
this sieve argument of ours. Thus,
if any in the family has a zero that violates (1.4),
then $\beta_T$ should not exist.
\medskip
Numerous extensions of the Deuring--Heilbronn
phenomenon and Siegel's theorem
have been considered, but the present work of ours
appears to be the
first concerning Linnik's.  Among the former is D.M. Goldfeld's
fundamental work [2]. He related 
exceptional zeros with the vanishing of the central values
of $L$-functions attached to elliptic 
curves and thus the Birch--Swinnerton-Dyer conjecture
came into the scene. An appropriate case of the conjecture
was confirmed later by B. Gross and D. Zagier (1986), 
which thus settled Gauss' conjecture mentioned 
above in an effective fashion.  However,
despite its amazing depth, their work
does not yield any  significant
consequence in the theory of 
the  distribution of primes, since it does not eliminate the possibility
of the existence of exceptional zeros and
the implied zero-free  region of 
Dirichlet $L$-functions is unfortunately
too weak  to be applied. 
\bigskip
\noindent
{\bf 2. Basic tools.}
Hereafter until the end of Section 4,
we shall assume that the $T$-exceptional 
character exists. We put
$$
\beta_T=1-\delta,\eqno(2.1)
$$
and introduce the multiplicative function
$$
f(n)=\sum_{d|n}\chi_T(d)d^{-\delta},\eqno(2.2)
$$
which is always positive. We consider the 
$\Lambda^2$-sieve situation
$$
\sum_{n<N}f(n)\Big(\sum_{d|n}\lambda_d\Big)^2,
\eqno(2.3)
$$
where $\lambda_1=1$ and $\lambda_d=0$ for
$d\ge z$, with a large parameter $z$.
The optimal $\{\lambda_d\}$ are given by
$$
\lambda_d=\mu(d)F_d{G_d(z/d)\over G_1(z)}
\eqno(2.4)
$$
where $\mu$ is the M\"obius function,
$$
F_d=\prod_{p|d}F_p,\quad
F_p=\Big(1-{1\over p}\Big)^{-1}
\Big(1-{\chi_T(p)\over p^{1+\delta}}\Big)^{-1},\eqno(2.5)
$$
and
$$
G_d(x)=\sum_{\scr{d<x}\atop\scr{(d,r)=1}}
{\mu^2(r)\over g(r)},\quad g(r)=\prod_{p|r}(F_p-1)^{-1},\eqno(2.6)
$$
with $p$ denoting a generic prime. 
The $\{\lambda_d\}$
defined by $(2.4)$ lead us to the 
multiplicative function
$$
\Phi_r(n)={\mu((r,n))g((r,n))},\eqno(2.7)
$$
where $(r,n)$ denotes the greatest common divisor of $r$ and $n$;
in fact, we have
$$
\sum_{d|n}\lambda_d ={1\over G_1(z)}\sum_{r<z}
{\mu^2(r)\over g(r)}\Phi_r(n).\eqno(2.8)
$$
The functions $\{\Phi_r\}$ have a 
quasi-orthogonality like Selberg's quasi-characters 
as is exhibited in
\medskip
\noindent
{\bf Lemma 1.}  {\it We have, for any 
complex numbers $\{a_n\}$ and for arbitrary $M, N,z$
with $N\ll M$, $z\ge T$,
$$
\eqalignno{
&\sum_{r<z}{\mu^2(r)\over g(r)}
\Big|\sum_{M\le n<M+N}\Phi_r(n)
f(n)a_n\Big|^2\cr
\le&\Big({\cal F}N+O(z^5M^{2/3})
\Big)\sum_{M\le n<M+N}f(n)|a_n|^2, &(2.9)
}
$$
where ${\cal F}=L(1+\delta,\chi_T)$.
}
\medskip
\noindent
{\it Proof\/}.  See [12, p.\ 33 and p.\ 187];
in fact a stronger result is given there, but this simplified version
suffices for our present purpose. 
The argument is, in essence, an application 
of the duality principle 
$\Vert{\cal D}\Vert=\Vert{\cal D}^*\Vert$ 
concerning an arbitrary bounded linear 
operator $\cal D$ and 
its adjoint ${\cal D}^*$ over a Hilbert space.
\medskip
\noindent
{\bf Lemma 2.} {\it 
We have, for any $z\ge T$,
$$
\delta^{-1}{\cal F}\ll G_1(z).\eqno(2.10)
$$
}
\par
\noindent
{\it Proof\/}. See [12,  p.\ 187]. The quantity
$1/G_1(z)$ represents the sieve effect of the optimal
$\{\lambda_d\}$.
\medskip
The following assertion is essential for our
purpose but is independent of the above:
\medskip
\noindent
{\bf Lemma 3.} {\it Let $v$ be a 
large positive parameter
and $\vartheta>0$
an arbitrary constant. We put, 
with an integer $l\ge1$,
$$
\Xi_d^{(l)}={1\over l!}(\vartheta\log v)^{-l}
\sum_{j=0}^l(-1)^{l-j}{l\choose j}\xi_d^{(j,l)},
\eqno(2.11)
$$
where
$$
\xi_d^{(j,l)}=\cases{\mu(d)
\Big(\displaystyle{\log {v^{1+j\vartheta}\over
d}\Big)^l}&$d<v^{1+j\vartheta}$,\cr\cr\hfil 0 
&$d\ge v^{1+j\vartheta}$.}\eqno(2.12)
$$
Then we have
$$
\Xi_d^{(l)}=\mu(d),\quad d<v\eqno(2.13)
$$
as well as
$$
\sum_{n=1}^\infty d_l(n)\Big(\sum_{d|n}
\Xi_d^{(l)}\Big)^2 n^{-\omega}\ll 1,
\eqno(2.14)
$$
whenever $ \omega\ge 1+1/\log v$; here $d_l(n)$ is as usual the
number of ways of representing $n$ as a product of $l$ 
integral factors.
}
\medskip
\noindent
{\it Proof\/}. This is Theorem 4 on [12, p.\ 24].

\bigskip
\noindent
{\bf 3. The hypothesis.}
We then deal with a multiplicative function 
$h$ that is to be convolved with $f$. 
We shall impose
conditions upon $h$. For this sake we need to
introduce the following notations:
$$
\eqalign{
F(s,h)&=\prod_p F_p(s,h),\quad
F_p(s,h)=\sum_{l=0}^\infty f(p^l)h(p^l) p^{-ls},\cr
H(s)&=\prod_p H_p(s),\quad
H_p(s)=\sum_{l=0}^\infty{h(p^l)p^{-ls}},\cr
H(s,\chi_T)&=\prod_p H_p(s,\chi_T),\quad
H_p(s,\chi_T)=\sum_{l=0}^\infty{\chi_T(p^l)h(p^l)p^{-ls}},
}\eqno(3.1)
$$
where $\Re s$ is assumed to be 
large initially. Then the conditions are:
\bigskip
\item{(A)} {\it There is an $L\in\B{N}$ such that $h(n)\ll d_L(n)$
for all $n\in\B{N}$};
\smallskip
\noindent
\item{(B)} {\it There exist constants 
${1\over2}<\eta<1$  and $\gamma>0$ such that 
the functions $H(s)$ and
$H(s,\chi_T)$ are regular and $\ll (T|s|)^\gamma$
for $\Re s>\eta$};
\smallskip
\item{(C)} {\it For any prime $p$, the 
functions $H_p(s)$, $H_p(s,\chi_T)$, and $F_p(s,h)$
do not vanish for $\Re s>\eta$};
\bigskip
\noindent
Since (C) is fulfilled by all sufficiently large primes, 
one may dispense with it. 
\medskip
\noindent
{\bf Lemma 4.} {\it Provided {\rm (A), (C)},
$\mu^2(r)=1$, and $\xi_d\ll\mu^2(d)$,
we have, for $\Re s>1$,
$$
\sum_{n=1}^\infty f(n)h(n)\Phi_r(n)
\Big(\sum_{d|n}\xi_d\Big)n^{-s}
=F(s,h)M_r(s,h;\xi),\eqno(3.2)
$$
where
$$
M_r(s,h;\xi)=g(r)\sum_{d=1}^\infty 
\xi_d\mu((r,d))
\prod_{p|d}(1-F_p(s,h)^{-1})
\prod_{\scr{p\nmid d}\atop\scr{p|r}}
(F_p(s,h)^{-1}F_p-1).
\eqno(3.3)
$$
}
\par
\noindent
{\it Proof\/}. This is a simple 
modification of Lemma 5 
on [12, p.\ 35], and the proof is analogous. 
\medskip
\noindent
{\bf Lemma 5.} {\it If $h$ in $(3.2)$ is replaced 
by the constant function $1$, then we have
$$
\sum_{r<v^{1+l\vartheta}}{\mu^2(r)\over g(r)}
\{M_r(1,1;\Xi^{(l)})\}^2\ll 1/(\c{F}\log v),\eqno(3.4)
$$
provided $\log v\gg\log T$.
\/}
\medskip
\noindent
{\it Proof\/}. This is Lemma 6 on [12, p.\ 36].

\bigskip
\noindent
{\bf 4. Main theorem.} We proceed to the proof of
our main assertion which is given at the end of this section.
We shall assume (A)--(C). 
\smallskip
First we note that for $\Re s>\eta$
$$
F(s,h)
=H(s)H(s+\delta,\chi_T)W_0(s)W_1(s),\eqno(4.1)
$$
where
$$
W_0(s)=\prod_{p <P_0}
{F_p(s,h)\over H_p(s)H_p(s+\delta,\chi_T)},\quad
W_1(s)=\prod_{p \ge P_0}
{F_p(s,h)\over H_p(s)H_p(s+\delta,\chi_T)},\eqno(4.2)
$$
with any sufficiently large constant $P_0$. 
The factor $W_1(s)$ is obviously regular and 
bounded for $\Re s>{1\over 2}$, and $W_0(s)$
is dealt with (C). Thus by (B) the function $F(s,h)$ is
regular for $\Re s>\eta$.
We set, in the previous sections,
$$
z=T,\; v=T^A,\; \vartheta=1,\; \xi_d=\Xi_d^{(l)},\eqno(4.3)
$$
with a sufficiently large $A$; 
actually $A=20$ suffices. The constant $l$ is to be
fixed later. Then, as it follows from $(3.3)$,
the (C) implies that $M_r(s,h;\Xi^{(l)})$ is 
regular and $\ll T^B$ for $\Re s>\eta$
with an effective constant $B$.
\smallskip
With this, let us consider the expression
$$
X_r={1\over 2\pi i}\int_{2-i\infty}^{2+i\infty}
F(\rho+w,h)M_r(\rho+w,h;\Xi^{(l)})
\Gamma(w)V^wdw.\eqno(4.4)
$$
Here $H(\rho)=0$, $\rho=\nu+i\tau$, with
$$
{1\over2}(\eta+1)\le\nu\le1,\; |\tau|<T,\eqno(4.5)
$$
and $V=T^C$ with a sufficiently large
$C$.  Note that we have
$$
F(\rho,h)=0,\eqno(4.6)
$$
because of $(4.1)$. We have, by $(2.13)$ and $(3.2)$, 
$$
X_r=e^{-1/V}+\sum_{n\ge v}\Phi_r(n)f(n)h(n)
\Big(\sum_{d|n}\Xi_d^{(l)}\Big)n^{-\rho}e^{-n/V}.
\eqno(4.7)
$$
We shift the line of integration in $(4.4)$ to 
$\Re w=-{2\over3}(\nu-\eta)$ 
and by $(4.6)$ get immediately
$$
X_r\ll V^{-(\nu-\eta)/2},\eqno(4.8)
$$
which implies that
$$
{1\over4}\le\Big|\sum_{v\le n\le  
V^2}\Phi_r(n)f(n)h(n)
\Big(\sum_{d|n}\Xi_d^{(l)}\Big)n^{-\rho}
e^{-n/V}\Big|^2.\eqno(4.9)
$$
We multiply both sides by the factor 
$\mu^2(r)/g(r)$ and sum over $r<z$, getting
$$
\eqalignno{
&{G_1(z)\over\log T}\ll \sum_{\scr{N=2^m}
\atop\scr{v\le N<V^2}}
\sum_{r<z}{\mu^2(r)\over g(r)}\cr
&\times\Big|\sum_{N\le n<2N}\Phi_r(n)f(n)h(n)
\Big(\sum_{d|n}\Xi_d^{(l)}\Big)n^{-\rho}
e^{-n/V}\Big|^2,&(4.10)
}
$$
where $m\in\B{N}$.
By virtue of Lemma 1, we have
$$
{G_1(z)\over\log T}
\ll {\cal F}V^{2(1-\nu)}
\sum_{n=1}^\infty f(n)|h(n)|^2
\Big(\sum_{d|n}\Xi_d^{(l)}
\Big)^2n^{-\omega_0},\eqno(4.11)
$$
with $\omega_0=1+(\log T)^{-1}$.
Further,  by $(2.14)$ with $l=2L^2$ we find that
$$
{G_1(z)\over \log T}\ll{\cal F}V^{2(1-\nu)}.\eqno(4.12)
$$
In view of Lemma 2, we have proved 
\medskip
\noindent
{\bf Theorem 1.} {\it We assume the existence of the $T$-exceptional
zero in the sense $(1.1)$. 
Then, under the hypothesis {\rm (A)--(C)}, all
the zeros of the function $H(s)$ in the region
$(4.5)$ satisfy $(1.4)$, provided
the constant $a_0$ is adjusted appropriately. \/}
\medskip
\noindent
A few remarks are in order: 
\smallskip
\par
\noindent
\item{(1)} We are able to
include the situation where $H(s)$ has a simple pole at $s=1$,
which occurs, for instance, 
if we consider $h\equiv 1$, i.e., $H(s)=\zeta(s)$. This violates
(B). Nevertheless, the above argument works well, 
since we have Lemma 5, although we skip details. Thus,
our theorem is applicable to Dedekind zeta and Hecke $L$-functions
of algebraic number fields as well. The case where $H(s)$ has
multiple poles at $s=1$ can also be included by an appropriate
modification of Lemma 5.
\item{(2)} Hecke $L$-functions associated with holomorphic cusp
forms on the hyperbolic upper half plane or more general $L$-functions
of similar nature can be included in Theorem 1, provided they
admit the twist by the real character $\chi_T$. 
On the other hand the situation with
Maass forms is an open question, as the condition (A) is then
hard to confirm. We surmise that there should be an appropriate
modification of Lemma 4.
\item{(3)}  It should be worth 
remarking that according to
A. Ogg [18] there are Dirichlet series with
Euler products which satisfy both (A) and (C) and
vanish at $s=1$. Thus, if we are ever able to apply
Theorem 1 to any of his functions, then the existence of
the exceptional zeros will be eliminated once and for all. 
However, what really matters is the confirmation of
the regularity condition (B), which does not seem
feasible.
We wonder nevertheless if it is absurd to try to find a function 
which satisfies (A)--(C) and vanishes at $s=1$. 
\item{(4)} It is possible to extend the notion of exceptional zeros to
those zeros lying close to the line $\Re s=1$
of any $L$-function which is not necessarily in Dirichlet's
family. In fact, one may extend Lemmas 1--4 for this purpose.
This line of consideration yielded even a novel 
way [12, Section 4.1] to discuss zero-free regions of
$\zeta(s)$; in fact, it gave later the
assertion [14] that appears to be beyond the reach of the
convexity argument of E. Borel and C. Carath\'eodory.
\bigskip
\noindent
{\bf 5. An extension of the Brun-Titchmarsh Theorem.}
In the rest we shall give a version of large sieve
extensions of the Brun--Titchmarsh theorem. This is in fact a
rework of our old file left unpublished since early 1980's 
which we originally intended to include into our 
lecture notes [12].
We publish it here, as it might have now
some fresh interest in the light of J. Maynard's recent work [8]
as well as what we have developed above.
\medskip
Let 
$$
\pi(x;k,\ell)=\sum_{\scr{p\le x}\atop\scr{p\equiv \ell\bmod k}}
1,\quad (k,\ell)=1. \eqno(5.1)
$$
Then Theorem 13 on [12, p.\ 140] asserts, among other things, that
we have, uniformly for $kQ^2\le x^{9/20-\varepsilon}$,
$$
\sum_{\scr{q\le Q}\atop\scr{(q,k)=1}}
\mathop{{\sum}^*}_{\chi\bmod q}
\bigg|\sum_{\scr{p\le x}\atop\scr{p\equiv\ell\bmod k}}
\chi(p)\bigg|^2\le{(2+o(1))x\over\varphi(k)
\log(x/(kQ^2)^{3/8})}\pi(x;k,\ell),
\eqno(5.2)
$$
provided $x$ is larger than a constant which is effectively
computable for each sufficiently small $\varepsilon>0$,
where $\varphi$ is the Euler totient function. 
In particular we have
$$
\pi(x;k,\ell)\le{(2+o(1))x\over\varphi(k)\log(x/k^{3/8})},
\quad k\le x^{9/20-\varepsilon},\eqno(5.3)
$$
which surpasses partly the famed bound
$$
\pi(x+y;k,\ell)-\pi(x;k,\ell)
\le {2y\over\varphi(k)\log(y/k)},\quad k<y,\eqno(5.4)
$$
due to H.L. Montgomery and R.C. Vaughan [9].
In contrast to this, Maynard [8] asserts in essence that
$$
\pi(x;k,\ell)\le{2x\over\varphi(k)\log x},\quad k\le x^{1/8},
\eqno(5.5)
$$
provided $x$ is larger than an effectively 
computable constant. He gives also a lower bound, though
we skip it in order to make our presentation simple; 
for the same reason, we also skip mentioning 
former improvements upon $(5.4)$ other than $(5.3)$.
The bound $(5.5)$ has been known as a kind of folklore 
among specialists, but with a much less precision
about the range of moduli. 
\medskip
We shall refine $(5.2)$ by
\medskip
\noindent
{\bf Theorem 2.} {\it There exists an effectively computable
constant $\Omega$ such that we have, uniformly for
$kQ\le x^\Omega$,
$$
\sum_{\scr{q\le Q}\atop\scr{(q,k)=1}}
\mathop{{\sum}^*}_{\chi\bmod q}
\bigg|\sum_{\scr{p\le x}\atop\scr{p\equiv\ell\bmod k}}
\chi(p)\bigg|^2\le{2x\over\varphi(k)\log x}\pi(x;k,\ell).
\eqno(5.6)
$$
}
\par
\noindent
{\it Proof\/}. This assertion is in fact a simple 
consequence of our version on [12, p.\ 185] of the
Linnik--Fogels--Gallagher prime number theorem mentioned above.
We put
$$
\psi(x,\chi)=\sum_{n\le x}\chi(n)\Lambda(n),\eqno(5.7)
$$
with the von Mangold function $\Lambda$;
and let
$$
\tilde{\psi}(x,\chi)=\cases{\psi(x,\chi)-x& if $\chi$ is principal,
\cr \psi(x,\chi)+x^{\beta_T}/\beta_T& if $\chi=\chi_T$,
\cr\psi(x,\chi)& otherwise.}\eqno(5.8)
$$
Then, Theorem 17, loc.cit., asserts that there exist effectively
computable absolute constants $a_1, a_2, a_3>0$ such that
provided $T^{a_3}\le x\le\exp((\log T)^2)$
$$
\sum_{q\le T}\mathop{{\sum}^*}_{\chi\bmod q}
|\tilde{\psi}(x,\chi)|\le a_1 x\Delta_T
\exp\big(-a_2\log x/\log T\big),\eqno(5.9)
$$
where
$$
\Delta_T=\cases{\hfil 1& if $\chi_T$ does not exist,\cr
(1-\beta_T)\log T & if $\chi_T$ exists.}\eqno(5.10)
$$
Then we note that
$$
\eqalignno{
\sum_{\scr{p\le x}\atop\scr{p\equiv\ell\bmod k}}
\chi(n)\Lambda(n)&={1\over\varphi(k)}\sum_{\xi\bmod k}
\bar\xi(\ell)\psi(x,\xi\chi)\cr
&={1\over\varphi(k)}\sum_{\xi\bmod k}\bar{\xi}(\ell)
\psi(x,\xi^\sharp\chi)+O(\nu(k)\log x),&(5.11)
}
$$
where $\xi^\sharp$ is the primitive character inducing 
the Dirichlet character $\xi$, and
$\nu(k)$ the number of distinct prime factors of $k$.
Here $\xi^\sharp\chi$ stands for a unique primitive
character whose conductor is not larger than
$kQ$. We have thus
$$
\eqalignno{
&\sum_{\scr{q\le Q}\atop\scr{(q,k)=1}}
\mathop{{\sum}^*}_{\chi\bmod q}
\bigg|\sum_{\scr{p\le x}\atop\scr{p\equiv\ell\bmod k}}
\chi(n)\Lambda(n)\bigg|\cr
&\le{1\over\varphi(k)}\left\{x+x^{\beta_0}/\beta_0
+E(x,kQ)\right\}+O(\nu(k)Q^2\log x),&(5.12)
}
$$
where $\beta_0$ is the $kQ$-exceptional zero if exists;
and $E(x,kQ)$ is the left side of $(5.9)$ for
$T=kQ$. If $\beta_0$ exists, then $(5.9)$ implies that
$$
\eqalignno{
x^{-1}\left(x^{\beta_0}/\beta_0+E(x,kQ)\right)\le&\,
\exp\big(-\Delta_T\log x/\log T\big)/
(1-\Delta_T/\log T)\cr
+&\,a_1\Delta_T\exp\big(-a_2\log x/\log T\big),&(5.13)
}
$$
provided $T^{a_3}\le x\le((\log T)^2)$. The right side is
$$
\eqalignno{
&\le\exp\big(-a_3\Delta_T)/(1-\Delta_T/\log T)
+a_1\Delta_T\exp(-a_2a_3)\cr
&\le\exp(-a_3\Delta_T)
+\Delta_T\big(1/(2\log T)+a_1\exp(-a_2a_3)\big)\cr
&<1-{1\over2}a_3\Delta_T,&(5.14)
}
$$
since we may assume that $a_3\Delta_T$ is small 
while $a_3$ is large. Hence we have proved that if $\beta_0$ exists,
then
$$
\sum_{\scr{q\le Q}\atop\scr{(q,k)=1}}
\mathop{{\sum}^*}_{\chi\bmod q}
\bigg|\sum_{\scr{p\le x}\atop\scr{p\equiv\ell\bmod k}}
\chi(n)\Lambda(n)\bigg|\le 2{x\over\varphi(k)}\Big
(1-{1\over5}a_3\Delta_T\Big),\eqno(5.15)
$$
provided $a_3\Delta_T$ is small and 
$x\ge T^c=(kQ)^c$ with an effective
absolute constant $c>0$.
The case where the $kQ$-exceptional zero does not exist
is analogous; in fact, simpler. The rest
of the proof may be skipped, as it is a routine application
of integration by parts.
\medskip
Obviously $(5.6)$ contains $(5.5)$ but for $k\le x^\Omega$. It 
remains thus to find a good lower bound 
for $\Omega$. We are certain that
Maynard's argument will extend to 
the direction indicated by $(5.6)$ and yield $(5.5)$ as a particular 
instance, since the basic structure of 
his argument is essentially the same as ours that is
developed in [11][12], although the intricate
part of [8] corresponding to the numerical
precision should be overhauled accordingly. 
Further, we add that it is possible
to prove a short interval version of $(5.6)$.
\medskip
If the coefficient 2 in $(5.2)$--$(5.6)$ 
is replace by any smaller effective constant, 
then the exceptional zeros should not exist; a proof can be found
in [12, Section 4.3]. This means that the coefficient 2 will then be
essentially halved immediately, which shows well a tantalising
nature of the problem of improving the Brun--Titchmarsh theorem.
\medskip
Here some more comments are in order:
The bounds $(5.2)$--$(5.4)$
are sieve results; 
that is, they are proved using mainly sieve arguments,
without the zero-density theory 
or the Linnik phenomenon. The proof in [12]
of the assertion (5.2) depends on Iwaniec's work [4] on the
bilinear structure in the error term of the combinatorial linear sieve;
an alternative approach to his
result itself can be found in [12] (see also [3]).
Prior to [4], a bilinear structure in the error term of the $\Lambda^2$
sieve was observed in [10] and the first
improvement upon $(5.4)$ was 
achieved; see [5] as well. Later
the development [13] made it possible to 
prove $(5.2)$ via the $\Lambda^2$ sieve; see [17] for 
a further development. On the other
hand, the bound $(5.9)$ depends on our large sieve extension
of the $\Lambda^2$ sieve that is devised 
via the duality principle and the quasi-character
derived from optimal $\Lambda^2$-weights, as is already
mentioned in the first section.
Finally we remark that Selberg [19, Section 22] improved
$(5.4)$ by replacing the denominator by
$\varphi(k)(\log (y/k)+2.8)$. This is, however, definitely weaker 
than $(5.3)$ in the range of variables indicated there. 
Apparently he did not notice the fact that $(5.3)$ was attainable 
with his $\Lambda^2$-sieve.
\bigskip
\noindent
\centerline{\bf References}
\medskip
\noindent
\item{[1]} E. Bombieri. {\it Le Grand Crible dans la
Th\'eorie Analytique des Nombres\/}.  (seconde 
\'edition revue et augment\'ee). Ast\'erisque {\bf18},
Soc.\ Math.\ France, Paris 1987/1974.
\item{[2]} D.M. Goldfeld. The class number of
quadratic fields and the conjectures
of Birch and Swinner\-ton-Dyer. Annali della Scuola
Norm.\ Sup.\ di Pisa Cl.\ Sci., (4) {\bf 3}, 
623--663 (1976).
\item{[3]} G. Greaves. {\it Sieves in Number Theory\/}. 
Springer-Verlag, Berlin 2001.
\item{[4]} H. Iwaniec. A new form of the error 
term in the linear sieve. Acta Arith.,  {\bf37} (1980), 307--320.
\item{[5]} ---. Sieve methods. In:
{\it Intern.\ Congress of Math.\ Proc., Helsinki 1978\/},  
Acad.\ Sci.\ Fennica, Helsinki 1980, pp.\ 357--364.
\item{[6]} M. Jutila. On Linnik's constant. Math.\ Scand., {\bf 41}
(1977), 45--62.
\item{[7]} A.A. Karatsuba. {\it Elements of Analytic Number 
Theory\/}. Nauka, Moscow 1975. (Russian)
\item{[8]} J. Maynard. On the Brun--Titchmarsh theorem.
arXiv:1201.1777 (math.NT).
\item{[9]} H.L. Montgomery and R.C. Vaughan. 
The large sieve. Mathematika, {\bf20} (1973), 119--134.
\item{[10]} Y. Motohashi. On some improvements of the 
Brun--Titchmarsh theorem. J. Math.\ Soc.\ Japan, 
{\bf 26} (1974), 306--323.
\item{[11]} ---. Primes in arithmetic progressions.
Invent.\ math., {\bf 44} (1978), 163--178.
\item{[12]} ---. {\it Sieve Methods and 
Prime Number Theory\/}. Lect.\ Notes in Math.\ Phys., {\bf72},
Tata IFR and Springer-Verlag, Bombay 1983. 
\item{[13]} ---. On the error term in the Selberg sieve.  
In: {\it Number Theory in Progress: A. Schinzel Festschrift\/}, 
Walter de Gruyter,  Berlin 1999, pp.\ 1053--1064.
\item{[14]} ---. An observation on the 
zero-free region of the Riemann
zeta-function. {\it A. S\'ark\"ozy Festschrift}, 
Periodica Math.\ Hungarica, 
{\bf 42} (2001), 117--122. 
\item{[15]} ---. {\it Analytic Number Theory\/}.\ I. 
{\it Distribution of Prime Numbers\/}. Asakura Books, 
Tokyo 2009; II. {\it Zeta Analysis\/}. ibid, 2011. (Japanese)
\item{[16]} ---. On some improvements of the 
Brun--Titchmarsh theorem.\ IV. arXiv:1201.3134v1 [math.NT].
\item{[17]} Y. Motohashi and J. Pintz.
A smoothed GPY sieve.
Bull.\ London Math.\ Soc., {\bf40} (2008), 
pp.\ 298--310. 
\item{[18]} A. Ogg. On a convolution of $L$-series. 
Invent.\ math., {\bf 7} (1969), 297--312.
\item{[19]} A. Selberg. Lectures on sieves. In:
{\it Collected Papers\/}, II,
Springer-Verlag, Berlin 1991, pp.\ 65--247.

\vskip 1cm
\noindent
Department of Mathematics,
\par
\noindent
Nihon University,
\par
\noindent
Surugadai, Tokyo 101-8308, JAPAN
\bye